\documentclass[11pt]{article}
\usepackage{graphicx}
\usepackage{amsmath}
\usepackage{amsfonts,color}
\usepackage{amssymb,url}

\def\i{\mathrm i}

\newtheorem{thm}{Theorem}[section]
\newtheorem{lem}[thm]{Lemma}
\newtheorem{prop}[thm]{Proposition}

\newtheorem{defn}[thm]{Definition}

\newtheorem{rem}[thm]{Remark}

\numberwithin{equation}{section}

\newcommand{\R}{\mathbb R}

\newcommand{\rr}{\mathbb R}
\newcommand{\E}{\mathbb E}

\newcommand{\Z}{\mathbb Z}
\newcommand{\N}{\mathbb N}
\newcommand{\eqfdd}{\stackrel{f.d.d.}{=}}
\newcommand{\ip}[2]{{\langle #1,#2\rangle}}

\renewcommand{\Im}{{\mathcal I}}
\newcommand{\1}{\mathbf 1}

\DeclareMathOperator{\M}{\mathbf{M}}
\DeclareMathOperator{\W}{\mathbf{W}}
\DeclareMathOperator{\e}{\mathbf{e}}

\DeclareMathOperator{\tr}{\textup{tr}}

\def\limfdd{\renewcommand{\arraystretch}{0.5}
\begin{array}[t]{c}
\stackrel{\rm f.d.d.}{\longrightarrow} \\
\end{array}\renewcommand{\arraystretch}{1}}

\newcommand{\longrightarrowD}{\stackrel{d}{\longrightarrow}}

\newcommand{\noi}{\noindent}

\begin{document}


\sloppy
\title{Asymptotic theory for regression models with fractional local to unity root errors}

\author{Kris De Brabanter  \  and \  Farzad  Sabzikar}
\date{\today \\  \small
\vskip.2cm
Iowa State University }

\maketitle

\begin{abstract}
This paper develops the asymptotic theory for parametric and nonparametric regression models when the errors have a fractional local to unity root (FLUR) model structure. FLUR models are stationary time series with semi-long range dependence property in the sense that their covariance function resembles that of a long memory model for moderate lags but eventually diminishes exponentially fast according to the presence of a decay factor governed by a noncentrality parameter. When this parameter is sample size dependent, the asymptotic normality for these regression models admit a wide range of stochastic processes with behavior that includes long, semi-long, and short memory processes.
\end{abstract}





\section{Introduction}

This paper discusses the asymptotic theory for parametric and nonparametric regression models when the errors follow an extension of the popular local to unity root (LUR) model.
A LUR model time series $\{ X(t)\}$ is generated by
\begin{equation}\label{eq:LTU}
X(t) = \rho_N X(t-1) + \zeta(t),\quad t=1,\ldots, N; X(0)=0,
\end{equation}
where $\rho_N= 1+\frac{c}{N}$ for $c\in\mathbb{R}$, $\{\zeta(j)\}_{j\in\mathbb{Z}}$ are i.i.d innovations, and $N$ is the sample size. We refer the interested reader to \cite{Bobkoski,ChanWei1,Phillips,Stock1,ElliottStock,CavanaghElliottStock,MoonPhillips,Rossi,CampbellYogo,JanssonMoreira,Mikusheva2007,Muller,MoonVelasco, Jeganathan1, Jeganathan2, Jeganathan3} for examples of the LUR model in economics. Using the backward shift operator $BX(t)=X(t-1)$, we rewrite \eqref{eq:LTU} as $(1-\rho_N B)X(t)=\zeta(t)$. This representation immediately leads to the following extension:
\begin{equation}\label{eq:FLTU}
(1-\rho_N B)^{d}X_{d,\rho_N}(t)=\zeta(t),
\end{equation}
where $d\in \mathbb{R}\setminus \{-1,-2,\ldots\}$. In what follows, we assume (i) $\rho_N= 1+c_N$ where $c_N\to 0$ as $N\to\infty$, (ii) $N(\rho_N-1)\to c_*\in [0,\infty]$.  Phillips \cite{Phillips2} called $c_N$ as  the noncentrality parameter. We call the FLUR model $\{X_{d,\rho_N}(t)\}$ strongly, weakly, and moderately FLUR model if $c_*=\infty, 0$, and $<\infty$ respectively. Using the binomial expansion, we write
\begin{equation}\label{eq:FLTUmoving}
X_{d,\rho_N}(t)=
\sum_{k=0}^{\infty} \rho^k_N\ \omega_{-d}(k)\zeta (t-k),\ t \in \mathbb{Z},
\end{equation}
where $\omega_{-d}(k)= \frac{\Gamma(k+d)}{\Gamma(k+1)\Gamma(d)}$.
On the other hand, Sabzikar and Surgailis \cite{SabzikarSurgailisinvariance} applied the idea of the tempered fractional difference operator in \cite{TFC} to introduce the tempered linear process (TLP) $\{ X_{d,\lambda}(j)  \}_{j\in\Z}$ with moving averages
\begin{equation}\label{semilongmemory}
X_{d,\lambda}(j)= \sum_{k=0}^{\infty} e^{-\lambda k}\ b_{d}(k) \zeta(j-k), \qquad j\in\Z,
\end{equation}
where $\{ \zeta(j)   \}_{j\in \Z}$ are i.i.d innovations, $b_d(k)$ regularly varying at infinity as $k^{d-1}$, viz.
\begin{equation} \label{bdk}
b_d (k)\ \sim \ \frac{c_d}{\Gamma(d)} \, k^{d-1}, \qquad k \to \infty,  \qquad c_d \ne 0,  \quad d \ne 0,
\end{equation}
where $ d \in\rr $ is a real number, $d \neq -1, -2, \ldots$, and $\lambda>0$ is the tempering parameter.
Therefore, we can interpret the FLUR model $X_{d,\rho_N}$ in \eqref{eq:FLTUmoving} as the tempered linear process $X_{d,\lambda_N}$.

The autocovariogram of FLUR and TLP resemble long range dependent series out to moderate lag lengths but eventually decays exponentially fast. Giraitis et al. \cite{giraitis} named this behavior as semi-long memory which is analogous to the semi-heavy tail property in \cite{B-N}. Giraitis et al. \cite{giraitissemilongRS} introduced semi-long memory FARIMA$(0,d,0)$, semi-long memory LARCH, and semi-long memory ARCH processes and used these models to investigate the power and robustness of the $R/S$ type tests under contiguous and semi-long memory alternatives. Dacorogna \cite{Dacorogna} and Granger and Ding \cite{GrangerDing} argued that the covariance function of some economic time series decay slowly at first but ultimately decay much faster, such as the magnitude of certain powers of financial returns.  As an example, \cite{Dacorogna} showed that the autocorrelation function of absolute (20 minute) returns on the USD-DEM exchange rate follows a hyperbolic decline for lags up to about 10 days, which is a characteristic property of exchange rate returns, but for long lags the autocorrelation decays more rapidly which indicates a characteristic of semi-long memory. For another empirical example, we refer to Sabzikar et al. \cite{ARTFIMAestimation}, where the semi-long memory model called ARTFIMA(0, 0.3, 0.025, 0), here $d=0.3, \lambda=0.025$, is
used to model the log returns for AMZN stock price from 1/3/2000 to 12/19/2017. The advantage of using the semi-long memory ARTFIMA is
the fact that we can capture aspects of the low frequency activity better than the ARFIMA time series in part of the long-
range dependence scenario. In the aforementioned AMZN example, the periodogram follows a power law at moderate
frequencies (ARFIMA(0, 0.3, 0)), but then levels off at low frequencies and the ARTFIMA(0,0.3,0.025,0) model was found
to be more appropriate in capturing this behavior.

Motivated by the aforementioned applications of the semi-long memory model, we develop the asymptotic theory for parametric and nonparametric regression models when the errors follow a FLUR model. In both parametric and nonparametric cases, we establish functional limit theorems for weighted sums of a FLUR model. More specifically, we obtain the limiting distributions of
\begin{equation}\label{nonparametricparametric_partial}
A_{1}(Nh, \lambda_N, d)\sum_{j=1}^{N} K\Big(\frac{Nx-j}{Nh}\Big) X_{d,\lambda_N}(j)\quad {\rm and}\quad A_{2}(N, \lambda_N, d)\sum_{j=-\infty}^{\infty} f\Big(\frac{j}{N}\Big) X_{d,\lambda_N}(j),
\end{equation}
where $K(\cdot)$ and $f(\cdot)$ are functions satisfying certain conditions which will be specified later, $0<x<1$, $h$ is the bandwidth, and $N$ is the sample size.
The assumption that $\lambda$ depends on $N$ is reflected in the scale factors $A_{1}$ and $A_{2}$ and also the limiting distributions of \eqref{nonparametricparametric_partial}. When $\lambda_*=\infty$,
\begin{equation}\label{A1A2strong}
A_{1}(Nh, \lambda_N, d)=\frac{\lambda_N^d}{\sqrt{Nh}} \quad {\rm and}\quad \ A_{2}(N, \lambda_N, d)=\frac{\lambda_N^d}{\sqrt{N}}
\end{equation}
for $d\in \rr\setminus\mathbb{N}_{-}$. For $(d,\lambda_*)\in (-\frac{1}{2}, \frac{1}{2})\times \{0\}$ and $(d,\lambda_*)\in (-\frac{1}{2},\infty)\times (0,\infty)$, the scale factors are
\begin{equation}\label{A1A2moderateweak}
A_{1}(Nh, \lambda_N, d)=\frac{1}{(Nh)^{d+\frac{1}{2}}} \quad {\rm and}\quad  \ A_{2}(N, \lambda_N, d)=\frac{1}{N^{d+\frac{1}{2}}}.
\end{equation}
For the strongly and weakly tempered cases, the limiting distributions of \eqref{nonparametricparametric_partial} can be written as the Wiener integrals $\int_\rr g(u) dB_{d}(u)$ where $B_{d}(u)$ is a  fractional Brownian motion (FBM) with the parameter $-\frac{1}{2}<d<\frac{1}{2}$, while for the moderate tempered case, the limiting distribution can be represented as $\int_\rr f(u) dB^{I\!I}_{d,\lambda_*}(u)$, where $B^{I\!I}_{d,\lambda_*}(u)$ is called tempered fractional Brownain motion of the second kind (TFBMII).
Next, we explain the main contributions of this work which will be discussed in Sections \ref{secnonparametric} and \ref{secparametric}. In Section \ref{secnonparametric}, we consider the nonparametric regression model in fixed-design
\begin{equation}\label{nonparametric_intro}
Y(j) = m\Big(\frac{j}{N}\Big) + X_{d,\lambda_N}(j) \qquad (j=1,\ldots, N),
\end{equation}
where $m(x)$ is the unknown regression function and $\{X_{d,\lambda_N}(j)\}_{j\in\mathbb{Z}}$ is a FLUR model. To estimate $m(x)$, we consider the Priestley-Chao kernel estimator~\cite{PriestleyChao}
\begin{equation}\label{estimator_intro}
\hat{m}(x) =\frac{1}{Nh} \sum_{j=1}^{N}K\Big(\frac{Nx-j}{Nh}\Big) Y(j).
\end{equation}
In Theorems \ref{flt} and \ref{asymptoticsNW}, we establish asymptotic normality, consistency, and convergence rate for the Priestley-Chao kernel regression estimator\cite{PriestleyChao} under different values of $\lambda_*\in[0,\infty]$. The results of Theorems \ref{flt} and \ref{asymptoticsNW} not only cover the asymptotic theory for nonparametric regression in
\cite{CsorgoMielniczuk3,CsorgoMielniczuk1,CsorgoMielniczuk2,Deo1,GuoKoul,Robinson}, but also extend to the case $d\geq 1/2$ when $\lambda_*\in (0,\infty)$.
When $\lambda_*\in (0,\infty)$, the asymptotic results includes continuous stochastic processes that have semi-long range dependence in the sense of \cite{giraitis}.

In Section \ref{secparametric}, we consider the parametric regression model
\begin{equation}\label{parametric_intro}
Y(j) = \mu\Big(\frac{j}{N}\Big) + X_{d,\lambda_N}(j) \qquad (j=1,\ldots, N),
\end{equation}
where $\mu(x)$ is a continuous polynomial regression function with unknown knots and error process $\{X_{d,\lambda_N}(j)\}_{j\in\mathbb{Z}}$. In our framework, we assume that the trend function $\mu(x)$ is given by a piecewise polynomial and we aim to make an inference on the unknown parameters of the polynomial function. Therefore, we shall make a distinction between our problem and the classical approximation or interpolation theory of splines or nonparametric
spline smoothing, see \cite{Boor,Wahba,Eubank,GreenSiverman}. In the aforementioned references, an unknown regression function was estimated by splines. However, our work is more related with (non)linear regression and change point problems. For a general overview on nonlinear regression, see~\cite{SeberWild} and~\cite{Ivanov}. For change point problems, we refer to~\cite{Brodsky} and~\cite{CsorgoHorvath}. We obtain a unified formula for the asymptotic distribution of the least squares estimator under different structures of the error process by using tempered fractional calculus, see Theorems~\ref{Theorem4} and~\ref{limitdistributiontheorem}.

The paper is organized as follows. In Section~\ref{secTFBMIIdefn}, we revise the definition and some important properties of TFBMII and its connection with tempered fractional calculus. In Section~\ref{secnonparametric} we establish the finite dimensional distribution of normalized partial sums of weighted tempered linear processes and use it to develop the asymptotic theory for the nonparametric regression model in \eqref{nonparametric_intro}. In Section~\ref{secparametric}, we establish asymptotic results for the least squares estimator of the unknown knots of a piecewise polynomial regression model in~\eqref{parametric_intro}. In Appendix, we first introduce tempered fractional calculus and then we develop the Wiener integrals with respect to TFBMII. All proofs can be found in Section~\ref{secproofs}.

In what follows, $C$ denotes generic constants which may be different at different locations. We write $\limfdd$ and $\eqfdd$  
for convergence and equality of distributions in the sense of finite-dimensional distributions respectively and $\longrightarrowD$ to show the convergence in distribution. Denote $\N_\pm := \{\pm 1, \pm 2, \dots \}, \,
\R_+ := (0, \infty), \,
(x)_\pm := \max (\pm x, 0),  x \in \R,  \, \int := \int_\R$ and let $L^p(\R) \, (p \ge 1)$ denote the Banach space of measurable functions $f: \R \to \R$ with finite norm $\|f\|_p = \big(\int |f(x)|^p dx \big)^{1/p} $.

\section{TFBMII: Definitions and main properties}\label{secTFBMIIdefn}
TFBMII is introduced and discussed in detail in \cite{SurgailisFarzadTFSMII}. For the sake of completeness, we give a short review of the definition and some essential properties of TFBMII that we will apply in the next sections.

Let $d>-\frac{1}{2}$, $\lambda>0$, and $\{B(t)\}_{t\in\rr}$ be Brownian motion with mean zero and variance $\sigma^2|t|$. A TFBMII can be defined as the Wiener integral
\begin{equation}\label{eq:TFBMII}
B^{I\!I}_{d,\lambda}(t):= \frac{1}{ \Gamma(d+1) } \int_{\rr} h_{d,\lambda}(t;y) dB(y),
\end{equation}
where the function $y \mapsto h_{d,\lambda}(t;y): \rr \to \rr $ is defined by
\begin{equation}\label{eq:integrand}
h_{d,\lambda}(t;y):=(t-y)_+^{d} e^{-\lambda (t-y)_+} - (-y)_+^{d} e^{-\lambda (-y)_+}+ \lambda \int_{0}^{t} (s-y)_{+}^{ d } e^{-\lambda(s-y)_{+}}\ ds.   
\end{equation}
Recall from\cite{MeerschaertsabzikarSPA} that the (positive and negative) tempered fractional integrals (TFI) and tempered fractional derivatives (TFD) of a function $f: \R \to \R $
are defined by
\begin{equation}\label{I}
{\mathbb I}^{\kappa,\lambda}_{\pm} f(y) :=\frac{1}{\Gamma(\kappa)}\int f(s)(y-s)_{\pm}^{\kappa-1}\ e^{-\lambda(y-s)_{\pm}}\ ds,
\quad \kappa >0
\end{equation}
and
\begin{equation}\label{D}
{\mathbb{D}}^{\kappa,\lambda}_{\pm}f(y) := {\lambda}^{\kappa}f(y)+
\frac{\kappa}{\Gamma(1-\kappa)}\int (f(y)-f(s))(y-s)_\pm^{-\kappa-1}
\,\ e^{-\lambda(y-s)_\pm}\ ds,  \quad 0< \kappa <1,
\end{equation}
respectively. The TFI in~\eqref{I} exists a.e. in $\R $ for each $f \in L^{p}(\R)$ and defines a bounded linear operator in $L^{p}(\R), p \ge 1$. The TFD in~\eqref{D} exists for any absolutely continuous function $f \in L^1 (\R)$ such that $f^\prime \in L^1(\R)$; moreover, it can be extended to the fractional Sobolev space $W^{\kappa,2}(\mathbb{R}):=\big\{ f\in L^{2}(\mathbb{R}): \int (\lambda^2 + \omega^2)^{\kappa} |\hat{f}(\omega)|^2\ \ d\omega<\infty \big\}$, see \cite{MeerschaertsabzikarSPA} for more properties of TFI and TFD. The function $h_{d,\lambda}(t;y)$ can be written as
\begin{eqnarray}
\Gamma (d+1)\ {\mathbb I}^{d,\lambda}_{-} {\bf 1}_{[0,t]}(y)
&=& h_{d,\lambda}(t;y),
\quad d>0, \label{I1}  \\
\Gamma (d+1)\ {\mathbb D}^{-d,\lambda}_{-} {\bf 1}_{[0,t]}(y)
&=&h_{d,\lambda}(t;y),
\quad -\frac{1}{2}< d < 0.
\end{eqnarray}
Therefore, we can represent TFBMII in \eqref{eq:TFBMII} as follows:
\begin{equation}\label{eq:connection with fractional calculus}
B^{I\!I}_{d,\lambda}(t) =
 \begin{cases}
\int   {\mathbb I}^{d,\lambda}_{-} {\bf 1}_{[0,t]}(y)\ dB(y), &  d>0, \\
\int   {\mathbb D}^{-d,\lambda}_{-} {\bf 1}_{[0,t]}(y)\ dB(y),  &  -\frac{1}{2}<d<0.
\end{cases}
\end{equation}
Using \eqref{eq:connection with fractional calculus}, Parseval's formula for stochastic integrals~\cite[Proposition 7.2.7]{SamorodnitskyTaqqu}) and the Fourier transform of TFI and TFD, we have the harmonizable representation of TFBMII
\begin{equation}\label{eq:har}
B^{I\!I}_{d,\lambda}(t)\eqfdd \frac{1}{\sqrt{2\pi}}\int_{\rr}\frac{e^{ i\omega t}-1}{\i\omega}(\lambda+ i\omega)^{-d}\ d\widehat{B}(\omega),
\end{equation}
where $\hat{B}$ is an even complex-valued Gaussian white noise, $\overline {\hat {B}(dx)} = {\hat B}(-dx) $ with  zero  mean and variance $\E |{\hat B}(dx)|^2 = dx$. The next proposition summarizes basic properties of $B^{I\!I}_{d,\lambda}(t)$. We refer the reader to\cite{SurgailisFarzadTFSMII} for more details.
\begin{prop}\label{lem:g squar integrable}
\smallskip
\noi (i) TFBMII $B^{I\!I}_{d,\lambda}$ in \eqref{eq:TFBMII} has stationary increments, such that
\begin{equation}\label{eq:scalingTFBMII}
\left\{ B^{I\!I}_{d,\lambda}(ct) \right\}_{t\in\rr}{\eqfdd}\left\{ c^{d+\frac{1}{2}} B^{I\!I}_{d,c\lambda }(t) \right\}_{t\in\rr}
\end{equation}
for any scale factor $c>0$ and is not a self-similar process.

\noi (ii) TFBMII $B^{I\!I}_{d,\lambda}$ in \eqref{eq:TFBMII} has a.s. continuous paths.

\noi (iii) For $d>0$, the covariance function of TFBMII $B^{I\!I}_{d,\lambda}$ is given by
\begin{equation}\label{covZ}
\E B^{I\!I}_{d,\lambda}(t)B^{I\!I}_{d,\lambda}(s) =
C(d,\lambda)\int_{0}^{t}\int_{0}^{s}
|u-v|^{ d-\frac{1}{2} }K_{ d-\frac{1}{2} }(\lambda|u-v|)dv\ du,
\end{equation}
where $C(d,\lambda)=\frac{2}{\sqrt{\pi}\Gamma(d)(2\lambda)^{d-\frac{1}{2}}}$, $d>0$ and $\lambda>0$. Here $K_{\nu}(x)$ is the modified Bessel function of the second kind~\cite[Chapter 9]{abramowitz}.
\end{prop}

\section{Nonparametric regression model in fixed equispaced design with FLUR errors}\label{secnonparametric}
In this section, we study the asymptotic theory for the unknown regression function by using the Priestley-Chao kernel regression estimator when the errors has a FLUR or tempered linear model structure. More specifically, we consider the nonparametric equispaced fixed design regression model
\begin{equation}\label{nonparametricregmodel}
Y(j) = m \Big(\frac{j}{N}\Big)+X_{d,\lambda_N}(j)\qquad (j=1,\ldots, N),
\end{equation}
where $\{ X_{d,\lambda_N}(j)\}_{j\in\Z}$ is satisfying \eqref{semilongmemory}-\eqref{bdk}. In addition, we assume that
\begin{eqnarray}
&&\sum_{k=0}^\infty k^j b_d(k) = 0,  \quad 0 \le j  \le \lfloor-d\rfloor,  \hskip.5cm
  -\infty  < d < 0,  \label{bdkneg} \\
&&\sum_{k=0}^\infty |b_d(k)| < \infty,  \quad c_0 := \sum_{k=0}^\infty b_d(k) \ne 0,  \qquad  d = 0 \label{bdk0}
\end{eqnarray}
Assumptions \eqref{bdkneg} and \eqref{bdk0} are necessary for the validity of the convergence results.

The main statistic considered in our framework is
\begin{equation}\label{eq:S definition}
S_{d,\lambda_N}(u):=\sum_{k=1}^{\lfloor Nu\rfloor }X_{d,\lambda_N}(k),  \qquad u\in [0,1].
\end{equation}
We assume the innovations $\{ \zeta(i) \}_{i\in\Z}$ in \eqref{semilongmemory} are i.i.d with zero mean and unit variance so that $ N^{ -\frac{1}{2} } \sum_{i=1}^{[Nt]} \zeta(i)\longrightarrowD B(t) $. Next, we state the required assumptions to obtain our main results in this section.

{\bf Assumption 1.} The tempering parameter $\lambda \equiv \lambda_N$ may depend on $N$ so that $\lambda_N=o(1)$ and
\begin{equation}\label{bandwidth and temering}
N\lambda \to \lambda_*  \in [0, \infty]
\end{equation}
as $N\to\infty$.

{\bf Assumption 2.} The bandwidth $h=h_N\to 0$ and $Nh\to\infty$ as $N\to\infty$.

{\bf Assumption 3.} Let $K$ be a symmetric density function with support on $[-1,1]$ with bounded first derivative $K'$.

{\bf Assumption 4.} $h\log(Nh)\to 0$ as $N\to\infty$

\smallskip

\begin{thm}\label{thm:the third question}
Let $\{ X_{d,\lambda_N}(j)\}_{j\in\Z}$ be tempered linear process satisfying \eqref{semilongmemory}-\eqref{bdk} and \eqref{bdkneg}-\eqref{bdk0}.Then under Assumptions 1-3 and for $0<x<1$, we have

\noi(a) if $K\in L^{2}(\rr)$, $\lambda_*=\infty$, and $d\in\rr\setminus \mathbb{N}_{-}$, then
\begin{equation}
\frac{\lambda_N^d}{\sqrt{Nh}}\sum_{j=1}^{N}K\Big(\frac{Nx-j}{Nh}\Big)X_{d,\lambda_N}(j) \limfdd  \int_{0}^{2} K^{'}(1-t)\ B(t) dt
\end{equation}
as $N\to\infty$.

\noi(b) if
\begin{equation}\label{A3zero}
K\in {\mathcal{A}_{d,0}} : = \Big\{ f\in L^{2}(\rr) : \int_\rr |\hat{f}(\omega)|^2 |\omega|^{-2d} d\omega <\infty \Big\},
\end{equation}
with $\lambda_*=0$, $-1/2<d<1/2$, then
\begin{equation}
\frac{1}{(Nh)^{d+1/2}}\sum_{j=1}^{N}K\Big(\frac{Nx-j}{Nh}\Big)X_{d,\lambda_N}(j) \limfdd  \Gamma^{-1}(d+1)\int_{0}^{2}K^{'}(1-t)\ B^{I\!I}_{d,0}(t) dt,
\end{equation}
as $N\to\infty$.

\noi(c) if
\begin{equation}\label{A3lambdaclass}
K\in {\mathcal{A}_{d,\lambda}} := \{f\in L^2(\R):\int_{\R} (\lambda^2+\omega^2)^{-d}|\hat f(\omega)|^2\,\ d\omega<\infty\},
\end{equation}
with $\lambda_*\in (0,\infty)$, and $d\in (0,\infty)$, then
\begin{equation}
\frac{1}{(Nh)^{d+1/2}}\sum_{j=1}^{N} K\Big(\frac{Nx-j}{Nh}\Big)X_{d,\lambda_N}(j) \limfdd  \Gamma^{-1}(d+1) \int_{0}^{2}K^{'}(1-t)\ B^{I\!I}_{d,\lambda_*}(t) dt
\end{equation}
as $N\to\infty$, where $B^{I\!I}_{d,\lambda_*}$ is a TFBMII.
\end{thm}

To estimate the regression function $m$, we consider the following kernel estimator~\cite{PriestleyChao}
\begin{equation}\label{estimator}
\hat{m}(x) =\frac{1}{Nh} \sum_{j=1}^{N}K\Big(\frac{Nx-j}{Nh}\Big)Y(j).
\end{equation}
Next, we establish asymptotic normality, consistency, and convergence rate for the estimator~\eqref{estimator}, see Theorems \ref{flt} and \ref{asymptoticsNW}.


\begin{thm}\label{flt}
Let $\{X_{d,\lambda_N}(j)\}_{j\in\Z}$ be tempered linear process satisfying \eqref{semilongmemory}-\eqref{bdk} and \eqref{bdkneg}-\eqref{bdk0}. Assume the tempering parameter $\lambda$, the bandwidth $h$ and kernel $K$ satisfy Assumption 1-3. Further, let $K \in L^1(\mathbb{R}) \cap L^2(\mathbb{R})$. Then, for fixed $x$, we have
\begin{itemize}
\item [(a)] If $\lambda_* =\infty$ and $d\in\rr\setminus{\mathbb{N}}_{-}$, then
\begin{equation*}
\frac{\lambda^d_N}{  \sqrt{Nh}  } \sum_{j=1}^N K\biggl(\frac{Nx-j}{Nh}\biggr)X_{d,\lambda_N}(j) \longrightarrowD N(0, \sigma^{2}_{d,\infty}),
\end{equation*}
where
\begin{equation}\label{sigmastrongly}
\sigma^{2}_{d,\infty} = \sigma^{2}\int_{-1}^{1} K^2(u) du
\end{equation}

\item [(b)] If $-1/2<d<1/2$ and $\lambda_* =0$, then
\begin{equation*}
\frac{1}{ (Nh)^{d+1/2}  } \sum_{j=1}^N K\biggl(\frac{Nx-j}{Nh}\biggr)X_{d,\lambda_N}(j)\longrightarrowD N(0, \sigma^{2}_{d,0}),
\end{equation*}
where $\sigma^{2}_{d,0} = {\rm Cov}\Big( \int_{\rr}K(u) dB_{d,0}(u),\int_{\rr}K(\nu) dB_{d,0}(\nu)\Big)$. If we restrict $0<d<1/2$, then
\begin{equation}\label{sigmaweakly}
\sigma^{2}_{d,0} = \sigma^{2}\int_{-1}^{1} \int_{-1}^{1} K(u) K(v) |u-v|^{2d-1} du\ dv.
\end{equation}

\item [(c)] If $d\in (0,\infty)$ and $0< \lambda_* <\infty$, then
\begin{equation*}
\frac{1}{ (Nh)^{d+1/2}  } \sum_{j=1}^N K\biggl(\frac{Nx-j}{Nh}\biggr)X_{d,\lambda_N}(j) \longrightarrowD N(0, \sigma^{2}_{d,\lambda_*}),
\end{equation*}
where
\begin{equation}\label{sigmamoderatly}
\sigma^{2}_{d,\lambda_*} = C(d,\lambda_*)\int_{-1}^{1} \int_{-1}^{1} K(u) K(v) |u-v|^{d-{1/2}} K_{ d-{1/2}} (\lambda_* |u-v| )du\ dv
\end{equation}
and $C(d,\lambda_*)$ is given by Proposition \ref{lem:g squar integrable} and $K_{\nu}(x)$ is the modified bessel function of the second kind.
\end{itemize}
\end{thm}

\begin{thm}\label{asymptoticsNW}
Let the assumptions of Theorem~\ref{flt} and Assumption 4 hold. For every $k\in\mathbb{N}$ and $0<x_1< \ldots < x_k<1$, we have
  \begin{itemize}
  \item[(a)] If $\lambda_* =\infty$ and $d\in\rr\setminus{\mathbb{N}}_{-}$, then
  \begin{equation*}
     \lambda^d_N \sqrt{Nh} \Big( \widehat{m}(x_1)-\E \widehat{m}(x_1), \ldots, \widehat{m}(x_k)-\E \widehat{m}(x_k) \Big) \longrightarrowD \sqrt{\sigma^2_{d,\infty}} \Big( N_1, \ldots, N_k\Big),
  \end{equation*}
where $N_1, \ldots, N_k$ are independent standard normal distributions and $\sigma^2_{d,\infty}$ is given by \eqref{sigmastrongly}.
  \item[(b)] If $-\frac{1}{2}<d<\frac{1}{2}$ and $\lambda_* =0$, then
  \begin{equation*}
     (Nh)^{1/2-d} \Big( \widehat{m}(x_1)-\E \widehat{m}(x_1), \ldots, \widehat{m}(x_k)-\E \widehat{m}(x_k) \Big) \longrightarrowD \sqrt{\sigma^2_{d,0}} \Big( N_1, \ldots, N_k\Big),
  \end{equation*}
where $N_1, \ldots, N_k$ are independent standard normal distributions and $\sigma^2_{d,0}$ is given by \eqref{sigmaweakly}.
  \item[(c)] If $d\in (0,\infty)$ and $0< \lambda_* <\infty$, then
  \begin{equation*}
     (Nh)^{1/2-d} \Big( \widehat{m}(x_1)-\E \widehat{m}(x_1), \ldots, \widehat{m}(x_k)-\E \widehat{m}(x_k) \Big) \longrightarrowD \sqrt{\sigma^2_{d,\lambda_*}} \Big( N_1, \ldots, N_k\Big),
  \end{equation*}
where $N_1, \ldots, N_k$ are independent standard normal distributions and $\sigma^2_{d,\lambda_*}$ is given by \eqref{sigmamoderatly}.
  \end{itemize}
\end{thm}

\begin{rem}\label{remextensionnonparametric}
{\rm When $\lambda_*=0$, the behaviour of $S_{d,\lambda_N}$ is typical for long range dependency. This fact explains why part (b) of Theorem~\ref{asymptoticsNW} is related with the asymptotic results for kernel regression function estimators obtained in~\cite[Proposition 1]{CsorgoMielniczuk1} and~\cite[Theorem 2]{Deo1} with long range dependence errors.}
\end{rem}

\begin{rem}\label{optimalrate}
\rm {
The tempering parameter $\lambda$ effects the optimal bandwidth and rate of convergence of the estimator. If the tempering parameter $\lambda$ does not depend on the sample size $N$, then applying the results in \cite{Hall} or \cite{BeranFeng2}, states that the optimal $L^{1}$ convergence rate for the regression estimator $\hat{m}(x)$ is of order $n^{-2/5}$ if a bandwidth of optimal order $n^{-1/5}$ is used. In fact,
according to~\cite[Eq. 3]{Ray}, the optimal bandwidth is
\begin{equation}\label{optimal_ARTFIMA}
h_{opt}= \Bigg\{\frac{(1-e^{-\lambda})^{-2d} \int_\rr K^2(x)\ dx}{\Big[\int_\rr {x^2} K(x) dx\Big]^2\int_0^1 [m^{''}(x)]^2\ dx}\Bigg\}^{1/5}N^{-1/5}.
\end{equation}
It would be interesting to consider the optimal bandwidth and the rate convergence when the tempering parameter $\lambda$ is sample size dependent. We conjuncture that depending on $\lambda_*\in [0,\infty]$ there are different optimal bandwidths. However, a rigorous simulation study would be needed and is beyond the scope of this paper.
Also, to the authors' knowledge, a method to estimate $\lambda_N$ does not exist. A possible lead could be found in~\cite{Debrabanter2018}. By using a kernel function $K$ such that $K(0)=0$, an asymptotically optimal bandwidth can be obtained by minimizing the residual sums of squares without any prior knowledge of the error process. We believe that this method could also be used to directly estimate the tempering parameter $\lambda_N$ from data.
}
\end{rem}

\section{Piecewise polynomial regression model with semi-long memory error}\label{secparametric}
In this section we investigate the asymptotic theory for piecewise polynomial and spline regression with partially unknown knots and errors having a FLUR or tempered linear model. We will obtain unified formulas
for the asymptotic distribution of least squares estimators of the unknown parameters based on tempered fractional calculus. More precisely, consider the parametric regression model
\begin{equation}\label{regmodel}
Y(j) = \mu \Big(\frac{j}{N}\Big)+X_{d,\lambda_N}(j)\qquad (j=1,\ldots, N),
\end{equation}
where $\mu(\cdot)$ is a continuous polynomial function such that
\begin{equation}\label{mfunctionparametric}
\mu(s)=\sum_{i=1}^{p}a_i f_{i}(s), \qquad s\in [0,1],
\end{equation}
and $\{ X_{d,\lambda}(j)\}$ is a tempered linear processes satisfying \eqref{semilongmemory}-\eqref{bdk} and \eqref{bdkneg}-\eqref{bdk0}. Here $a=(a_1,\ldots, a_p)$ denotes unknown regression coefficients and $f_1, \ldots, f_p$ are truncated power spline basis functions defined as $f_1(s)=1, f_2(s)=s, \ldots, f_q(s)=s^{q-1}, f_{q+1}(s)= (s - \eta)_+, \ldots, f_p(s)=(s - \eta)_{+}^{p-q}$ with $(s-\eta)_+ :={\bf 1}_{\{s-\eta>0\}}$ and $\eta$ is an unknown knot parameter. The regression function $\mu(s) =\sum_{i=1}^{p} a_i f_{i}(s)=\mu(s; \theta )$
depends on a $(p+1)$-dimensional parameter vector $\theta =(a^{T},\eta)^{T}\in \Theta=\rr^{p}\times (0,1)$. Since we want the model to be identifiable, we assume in addition $a_{i}\neq 0$ for at least one
$i\in \{q + 1, \ldots, p\}$.

Let $\hat{\theta}$ be the ordinary least squares estimator of $\theta$ that minimizes
\begin{equation*}
\beta(\tau) = \sum_{j=1}^{N} \Biggl[ Y(j)-\mu\biggl(\frac{j}{N}, \theta\biggr) \Biggr]^{2}
\end{equation*}
with respect to $\theta\in\Theta$. Let $\eta\in (0, 1)$ and define the $n \times p$ matrix
\begin{equation*}
\W_N = \W_N(\eta)=(w_{j,i} )_{j=1,\ldots,N; i=1,\ldots, p} = (\mathbf{w}_{N,1}, \ldots, \mathbf{w}_{N,p})
\end{equation*}
with $w_{j,i} = f_i \big(\frac{j}{N}\big) (1 \leq j \leq N; 1 \leq i \leq p)$, and column vectors denoted by $w_{N,i} (i = 1, \ldots, p)$. According to~\cite{Beran2}, $\W^{T}_{N}\W_{N}$ is invertible for large $N$ and the projection matrix onto the column space of $\W_N(\eta)$ is
\begin{equation*}
P_{\W_{N}} = P_{\W_N(\eta)} = {\W_N(\W_N^{T}\W_N)^{-1}\W_N^{T}}.
\end{equation*}
Hence, for given observations $X = (X_1, \ldots, X_n)^{T}$, $\hat{\eta}$ is obtained by minimizing $\beta = \| X - P_{ W_N( \tilde{\eta} ) } X\|_{2}$ with respect to $\tilde{\eta}$. Now, one can compute $ \hat{a} = ( \hat{a}_1, \ldots, \hat{a}_p )$  by projecting $X$ onto the column space of the design matrix $\W_N (\hat{\eta})$. For the estimated mean function we have
\begin{equation*}
\Big[
\mu\Big( \frac{1}{N}, \hat{\theta} \Big), \mu\Big( \frac{2}{N},\hat{\theta} \Big),
\ldots, \mu\Big( 1,\hat{\theta} \Big)\Big]^{T}
= P_{\W_N (\hat{\eta} ) }{\bf X} = P_{\W_{N}(\hat{\eta})} [ \mu_N(\theta) + {\bf e}_{N} ],
\end{equation*}
where $\mu_N (\theta) = [ \mu(\frac{j}{N}, \theta) ]_{j=1,\ldots, N}$ and ${\bf e}_N = ( X_{d,\lambda_N}(1), \ldots, X_{d,\lambda_N}(N) )^{T} $. In general, the partial derivatives $\partial_t \mu(t;\eta, a)$ and $\partial_\eta \mu(t;\eta, a)$ do not exists if $\eta=t$. However, one can use the left and right derivatives of $\mu(t; \eta, a)$ since they exist everywhere. We denote the left and right partial derivatives of $\mu$ with respect to $\theta_i$ by $\mu_{(i-)}$ and $\mu_{(i+)}$ respectively. Partial derivatives in the sense of an absolutely continuous function will be denoted by $\mu_{(i)}$. Note that $\mu_{(i+)} = \mu_{(i-)} = \mu_{(i)}$ almost everywhere. Defining the $N \times (p+1)$ matrix
\begin{equation*}
\M_{N+} = [\mu_{(i+)}(t/N)]_{t=1,\ldots,n;i=1,\ldots,p+1} \in \rr^{N\times (p+1)},
\end{equation*}
we have
\begin{equation*}
\lim_{N\to\infty} N^{-1}(\M_{N+}^T\M_{N+})_{jk}=\int_0^1 \mu_{(j)}(s,\theta) \mu_{(k)}(s,\theta)\ ds.
\end{equation*}
Similarly to $\W_{N}^T\W_{N} $, the matrix $\M_{N+}^T\M_{N} $ has full rank for large $N$ such that
\begin{equation}\label{biglambda}
{\bf \Lambda}=\lim_{N\to\infty} N\ (\M_{N+}^T\M_{N+})^{-1}
\end{equation}
is well defined. Theorem~\ref{Theorem1} shows that $\lambda_N^d \sqrt{N}(\hat{\theta}-\theta)$ and $N^{\frac{1}{2} - d}(\hat{\theta}-\theta)$ are asymptotically equivalent to $\lambda_N^d \sqrt{N}(\M_{N+}^T\M_{N+})^{-1} \M_{N+}^T\e_{N}$ and $N^{\frac{1}{2} - d}(\M_{N+}^T\M_{N+})^{-1} \M_{N+}^T\e_{N}$ for different values of $\lambda_\star$ respectively. Theorems~\ref{Theorem4} and~\ref{limitdistributiontheorem} establish the asymptotic distribution of the parameter vector $\widehat{\theta}$ for different values of $\lambda_\star$.
\begin{thm}\label{Theorem1}
Let $ \{ X_{d,\lambda_N}(t) \}_{t\in\Z}$ be a tempered linear process given by \eqref{semilongmemory}.
Then for any $\Delta>0$
\begin{itemize}
\item [(a)] If $\lambda_*=\infty$, $d>0$, and $\lambda_N= o(N^{-1/(2-2d)})$, then
\begin{equation}\label{eq:Theorem_3.1_Beran_strong}
\mathbb{P}\Big( \lambda_N^d \sqrt{N} \Big\|\hat{\theta}-\theta- (\M_{N+}^T\M_{N+})^{-1} \M_{N+}^T\e_{N}\Big\| >\Delta \Big)=o(1),
\end{equation}
as $N\to\infty$,
\item [(b)] If $\lambda_*\in [0, \infty)$ and $0< d <\frac{1}{2}$, then
\begin{equation}\label{eq:Theorem_3.1_Beran}
\mathbb{P}\Big( N^{\frac{1}{2} - d}\ \Big\|\hat{\theta}-\theta- (\M_{N+}^T\M_{N+})^{-1} \M_{N+}^T\e_{N}\Big\| >\Delta \Big)=o(1),
\end{equation}
as $N\to\infty$, where $\|\cdot\|$ denotes the Euclidean norm.
\end{itemize}
\end{thm}

\begin{thm}\label{Theorem4}
Let $ \{ X_{d,\lambda_N}(t) \}_{t\in\Z}$ be a tempered linear process given by \eqref{semilongmemory}.
\begin{itemize}
\item [(a)] If $d\in \mathbb{R}\setminus\mathbb{N}_{-}$ and $\lambda_*=\infty$, then
\begin{equation*}
\lambda_N^d \sqrt{N} \big( \hat{\theta}-\theta \big) \limfdd \Lambda \Xi,
\end{equation*}
where ${\bf \Lambda}$ is given by \eqref{biglambda}, $ \Xi = \Big[\int_{\rr} \mu_{(i+)}(s) dB(s)\Big]_{i=1,\ldots, p+1}$ is a random vector process, and $B(s)$ is a Brownian motion.

\item [(b)] If $0<d<1/2$ and $\lambda_*=0$, then
\begin{equation*}
N^{\frac{1}{2} - d }\big( \hat{\theta}-\theta \big) \limfdd \Lambda \Xi,
\end{equation*}
where ${\bf \Lambda}$ is given by \eqref{biglambda}, $ \Xi = \Big[\int_{\rr} \mu_{(i+)}(s)\ dB^{I\!I}_{d,0}(s)\Big]_{i=1,\ldots, p+1}$ is a random vector process, and $B^{I\!I}_{d,0}$ is a multiple of FBM.
\item [(c)] If $0<d<1/2$ and $\lambda_* \in (0,\infty)$, then
\begin{equation*}
N^{\frac{1}{2} - d }\big( \hat{\theta}-\theta \big) \limfdd \Lambda \Xi,
\end{equation*}
where ${\bf \Lambda}$ is given by \eqref{biglambda}, $ \Xi = \Big[\int_{\rr} \mu_{(i+)}(s)\ dB^{I\!I}_{d,\lambda_*}(s)\Big]_{i=1,\ldots, p+1}$ is a random vector process, and $B^{I\!I}_{d,\lambda_*}$ is TFBMII.
\end{itemize}
\end{thm}

The following theorem shows that the limit distribution in Theorem \ref{Theorem4} is Gaussian and also gives the closed form of the covariance matrix of the random vector $\Xi$.
\begin{thm}\label{limitdistributiontheorem}
Under the assumptions of Theorem~\ref{Theorem4}, the limit distribution of the random vector $\Lambda \Xi$ is Gaussian with zero mean and covariance $\Lambda \Sigma_0 \Lambda$. That is
\begin{equation*}
N^{\frac{1}{2} - d }\big( \hat{\theta}-\theta \big) \limfdd \Lambda \Xi \sim N(0, \Lambda \Sigma_0 \Lambda)
\end{equation*}
as $N\to\infty$. Moreover,
\begin{itemize}
\item [(a)] for $d\in \mathbb{R}\setminus\mathbb{N}_{-}$ and $\lambda_*=\infty$ and $0<d<\frac{1}{2}$, we have
\begin{equation}\label{SigmaBM}
\Sigma_{\infty}=\Bigg [  \int_\rr \int_\rr \mu_{(i+)}(t) \mu_{(k+)}(s) ds \ dt.   \Bigg]_{i,k=1,\ldots,p+1}.
\end{equation}

\item [(b)] for $0<d<1/2$ and $\lambda_*=0$ and $0<d<\frac{1}{2}$, we have
\begin{equation}\label{SigmaFBM}
\Sigma_{0}=\Bigg [\int_{\mathbb{R}} \Big( {\mathbb I}^{d,0}_{-} \mu_{(i+)} \Big)(s)
 \Big( {\mathbb I}^{d,0}_{-} \mu_{(k+)} \Big)(s)\ ds\Bigg]_{i,k=1,\ldots,p+1}.
\end{equation}
Moreover the covariance $\eqref{SigmaFBM}$ can be written as
\begin{equation}
\Sigma_0 = \frac{2}{\Gamma(d) \sqrt{\pi} (2\lambda)^{d-\frac{1}{2}}}\int_\rr \int_\rr \mu_{(i+)}(t) \mu_{(k+)}(s) |t-s|^{2d-1}  ds \ dt.
\end{equation}

\item [(c)] for $0<d<1/2$ and $\lambda_*\in(0,\infty)$ and $0<d<\frac{1}{2}$, we have
\begin{equation}\label{Sigma}
\Sigma_{\lambda_*}=\Bigg [\int_{\mathbb{R}} \Big( {\mathbb I}^{d,\lambda_*}_{-} \mu_{(i+)} \Big)(s)
 \Big( {\mathbb I}^{d,\lambda_*}_{-} \mu_{(k+)} \Big)(s)\ ds\Bigg]_{i,k=1,\ldots,p+1}.
\end{equation}
Moreover the covariance $\eqref{Sigma}$ can be written as
\begin{equation}
\Sigma_{\lambda_*} = \frac{2}{\Gamma(d) \sqrt{\pi} (2\lambda)^{d-\frac{1}{2}}}\int_\rr \int_\rr \mu_{(i+)}(t) \mu_{(k+)}(s) |t-s|^{d-\frac{1}{2}} K_{d-\frac{1}{2}}(\lambda_*|t-s|) ds \ dt.
\end{equation}
\end{itemize}
\end{thm}

\begin{rem}
{\rm \noi(a) Parts (a) and (b) of Theorem \ref{limitdistributiontheorem} are related to the result in \cite[Theorem 3.3]{Beran2} for the short and long memory cases.

\noi(b) Theorem \ref{limitdistributiontheorem} is also valid for $-\frac{1}{2} < d < 0$ and $\lambda_* \in [0,\infty)$. In this case,
\begin{equation}
\Sigma_{0}=\Bigg [\int_{\mathbb{R}} \Big( {\mathbb D}^{-d,\lambda_*}_{-} \mu_{(i+)} \Big)(s)
 \Big( {\mathbb D}^{-d,\lambda_*}_{-} \mu_{(k+)} \Big)(s)\ ds\Bigg]_{i,k=1,\ldots,p+1}.
\end{equation}
However, a closed form for $\Sigma_0$ for the case $d>0$ does not exist.

}
\end{rem}

\begin{rem}
{\rm All results of this section consider one unknown knot. However, similar results can be obtained for an arbitrary continuous piecewise polynomial function
\begin{equation*}
\mu(s) = \sum_{k=0}^{l} \sum_{j=1}^{p_k} a_{k,j} (s -\eta_k)_+^{b_{j,k} }
\end{equation*}
with $b_{j,k} < b_{j+1,k} $, knots $ 0 = \eta_0 < \eta_1 < \ldots < \eta_l < 1$ of which some (but not necessarily all) are unknown, and the condition $b_{j,k} \geq 1$ for $k\geq 1$
(needed for continuity).   }
\end{rem}

\appendix
\section{Proofs}\label{secproofs}
Before we prove the main results of the paper, we first state two technical lemmas upon which our results are based. Lemma \ref{lem1} and Lemma~\ref{thm:similar the third question}, play an important role in establishing the asymptotic results in Section~\ref{secparametric}. Next, we introduce some notations that will be used in Lemma~\ref{thm:similar the third question}.

For the function $f$ and $m\in\mathbb{N}\cup\{\infty\}$, we define the approximation
\begin{equation*}
f^{+}_{N,m}(y)=\sum_{j=0}^{m} f\Big(\frac{j}{N}\Big)1_{[\frac{j}{N},\frac{j+1}{N}]}(y), \qquad f^{-}_{N,m}=\sum_{j=-m}^{-1} f\Big(\frac{j}{N}\Big)1_{[\frac{j}{N},\frac{j+1}{N}]}(y),
\end{equation*}

\begin{equation*}
f^{+}_{N}=f^{+}_{N,\infty}, \qquad f^{-}_{N}=f^{+}_{N,\infty},\qquad f_{N}=f^{+}_{N}+f^{-}_{N}.
\end{equation*}

\begin{lem}
\label{lem1}
  Let $b_d(k)$ be defined as in~\eqref{bdk}. Then any $y \in \mathbb{R}$ we have
  \begin{equation}
  \biggl(\lambda_N + i\frac{y}{N}\biggr)^{d} \sum_{k=0}^{\infty} e^{- (\lambda_N + i\frac{y}{N})k}\, b_{d}(k) \sim 1
  \label{tauberian}
  \end{equation}
as $N\to\infty$.
\end{lem}

{\bf proof Lemma~\ref{lem1}}: For $d>0$, since $\sum_{k=0}^Nb_d(k) \sim 1/(d\Gamma(d)) N^d$ as $N\to \infty$ according to~\eqref{bdk} and $e^{-(\lambda_N+\frac{ix}{N})} \leq 1$, then according to the Tauberian theorem for power series~\cite[p. 447 Theorem 5]{Feller1971} we have
\begin{equation*}
\sum_{k=0}^{\infty} e^{- (\lambda_N + i\frac{y}{N})k}\, b_{d}(k) \sim (1-e^{- (\lambda_N + i\frac{y}{N})k})^{-d} \textrm{ as } N\to\infty
\end{equation*}
and consequently $\bigl(\lambda_N + i\frac{y}{N}\bigr)^{d}/(1-e^{- (\lambda_N + i\frac{y}{N})k})^{d}\sim 1$ for $\lambda_N \to 0$ and $N \to \infty$ as $N\to\infty$, proving~\eqref{tauberian}. For $-1<d<0$, define $\widetilde{b}_d(k) = \sum_{i=k}^\infty b_d(i) \sim -1/(d\Gamma(d))k^{\widetilde{d}-1}$ with $\widetilde{d}=d+1 \in (0,1)$. Next, we have that
\begin{equation}
 \widetilde{b}_d(0) = \sum_{i=0}^\infty b_d(i) = \sum_{i=0}^{k-1} b_d(i) + \sum_{i=k}^\infty b_d(i) = 0
 \label{switch}
\end{equation}
and therefore $\sum_{i=0}^{k-1} b_d(i) = -\sum_{i=k}^\infty b_d(i)$. Using summation by parts~\cite[p. 32, Eq. 2.5.8]{koul} yields
\begin{eqnarray*}
\lim_{s\to \infty} \sum_{j=0}^s e^{- (\lambda_N + i\frac{y}{N})j}\, b_{d}(j)&=& \lim_{s\to \infty} \bigl[e^{-(\lambda_N +i\frac{y}{N})s} \sum_{j=0}^s b_d(j)\\ && +  \sum_{j=0}^{s-1}e^{-(\lambda_N +i\frac{y}{N})s}-e^{-(\lambda_N +i\frac{y}{N})(s+1)} \sum_{i=0}^j b_d(i)\bigr] \\
&=& 0 + \{1- e^{- (\lambda_N + i\frac{y}{N})}\}\lim_{s\to \infty} \sum_{j=0}^{s-1}e^{-(\lambda_N +i\frac{y}{N})j} \sum_{i=0}^j b_d(i).
\end{eqnarray*}
By setting $j=k-1$ and using~\eqref{switch} we have
\begin{eqnarray*}
\lim_{s\to \infty} \sum_{j=0}^s e^{- (\lambda_N + i\frac{y}{N})j}\, b_{d}(j) &=& \{1- e^{- (\lambda_N + i\frac{y}{N})}\}\lim_{s\to \infty} \sum_{k=1}^{s-1}e^{-(\lambda_N +i\frac{y}{N})(k-1)} \sum_{i=0}^{k-1} b_d(i) \\
&=& -e^{(\lambda_N + i\frac{y}{N})} \{1- e^{- (\lambda_N + i\frac{y}{N})}\}\lim_{s\to \infty}\sum_{k=1}^{s-1}e^{-(\lambda_N +i\frac{y}{N})k}\sum_{i=k}^\infty b_d(i) \\
&=& -e^{(\lambda_N + i\frac{y}{N})} \{1- e^{- (\lambda_N + i\frac{y}{N})}\}\lim_{s\to \infty}\sum_{k=1}^{s-1}e^{-(\lambda_N +i\frac{y}{N})k}\widetilde{b}_d(k).
\end{eqnarray*}
Application of the Tauberian theorem for power series~\cite[p. 447 Theorem 5]{Feller1971} yields
\begin{equation*}
 \sum_{k=0}^{\infty} e^{- (\lambda_N + i\frac{y}{N})k}\, b_{d}(k) \sim  \{1- e^{- (\lambda_N + i\frac{y}{N})}\}^{1-\widetilde{d}} = \{1- e^{- (\lambda_N + i\frac{y}{N})}\}^{-d}
\end{equation*}
as $N\to \infty$, proving ~\eqref{tauberian}. In the general case $-j < d < -j +1$, $j = 1,2,\ldots $~\eqref{tauberian} follows similarly using summation by parts $j$ times. For $d=0$, it can be shown that the same result holds under an additional assumption on the sum of the $b_d(k)$'s~\cite{SabzikarSurgailisinvariance}. \hfill $\Box$

\begin{lem}\label{thm:similar the third question}
Let $\{ X_{d,\lambda_N}(j)\}_{j\in\Z}$ be tempered linear process given
by \eqref{semilongmemory}, $N\lambda_N\to\lambda_*\in (0,\infty)$ and $d>-1/2$. Let ${\mathcal{A}}_{d,\lambda_*}$ be the class of functions defined by \eqref{A3lambdaclass} and let
\begin{equation*}
{\rm{\bf Condition\ A}:}\ f,f^{\pm}_{N}\in{ {\mathcal{A}}_{d,\lambda_N}}, \|f^{\pm}_{N} - f^{\pm}_{N,m}\|_{ {\mathcal{A}}_{d,\lambda_N} }\to 0,\ {\rm as}\ m\to\infty,\ \|f-f_N\|_{ {\mathcal{A}}_{d,\lambda_*} }\to 0,\ {\rm as}\ N\to\infty
\end{equation*}
be satisfied, then
\begin{equation}
\frac{1}{N^{d+1/2}}\sum_{j=-\infty}^{\infty}f\Big(\frac{j}{N}\Big)X_{d,\lambda_N}(j) \limfdd  \int_{\rr}f(u)\ dB^{I\!I}_{d,\lambda_*}(u)
\end{equation}
as $N\to\infty$.
\end{lem}

\noindent{\bf proof Lemma~\ref{thm:similar the third question}:}
We first note that $\{X_{d,\lambda_N}(j)\}_{j\in\Z}$ can be written as
\begin{equation}\label{spectralreptempered}
X_{d,\lambda_N}(j) = \frac{1}{\sqrt{2\pi}} \int_{-\pi}^{\pi}e^{i\omega j}\sum_{k=0}^{\infty} e^{-i\omega k} e^{-\lambda_N k} b_{d}(k) \hat{B}(d\omega),
\end{equation}
where $\hat B(d \omega) $ is complex-valued Gaussian noise with $\E |\hat B(d \omega)|^2 = d \omega $, see \cite[Sections 4.6-4.7]{BrockwellDavisTSTM}.
Define
\begin{equation}\label{eq:UN and U defn}
U_{N}=\frac{1}{N^{ d+ \frac{1}{2}}} \sum_{j= -\infty}^{\infty} f\Big(\frac{j}{N}\Big) X_{d,\lambda_N}(j),
\qquad U=\int_{\rr}\ f(u)\ dB^{I\!I}_{d,\lambda_*}(u).
\end{equation}
The Wiener integral $U$ is well-defined, since $f\in{\mathcal{A}}_{d,\lambda}$. To show that the series $U_{N}$ is well-defined in the $L^{2}(\Omega)$, first apply the spectral representation of $\{X_{d,\lambda_N}(j)\}_{j\in\Z}$ given by \eqref{spectralreptempered}
\begin{equation}\label{sumwelldefined}
\begin{split}
\frac{1}{ {N}^{d+ \frac{1}{2}} } \sum_{j=0}^{m} f\Big(\frac{j}{N}\Big) X_{d,\lambda_N}(j)& = \frac{1}{ {N}^{d+ \frac{1}{2}} } \int_{-\pi}^{\pi} \Bigg[\sum_{j=0}^{m} \frac{1}{\sqrt{2\pi}} f\Big(\frac{j}{N}\Big) e^{ij\omega}\Bigg]\sum_{k=0}^{\infty} e^{ -(\lambda_N + i\omega )k } b_{d}(k)\ d\widehat{B}(\omega)\\
&=\frac{1}{ {N}^{d+ \frac{1}{2}}   } \int_{\rr}\Bigg[\sum_{j=0}^{m} \frac{1}{\sqrt{2\pi}} f\Big(\frac{j}{N}\Big) e^{\frac{i j y}{N}} \Bigg]{\bf 1}_{[-N\pi ,N\pi]}(y)\\
&\qquad\qquad\qquad\times \sum_{k=0}^{\infty} e^{ -(\lambda_N + i\frac{y}{N} )k } b_{d}(k) d\hat{B}\big(N^{-1}y\big)\\
&=\frac{1}{N^{d-\frac{1}{2}}}\int_{\rr} \Bigg[ \sum_{j=0}^{m} \frac{1}{\sqrt{2\pi}} f\Big(\frac{j}{N}\Big) \frac{e^{\frac{i(j+1)y}{N}}-e^{\frac{ijy}{N}}}{iy}
\Bigg]\\
&\qquad\qquad\qquad\times \frac{\frac{iy}{N}}{e^{\frac{iy}{N}}-1} 1_{[-N\pi ,N\pi]}(y) \sum_{k=0}^{\infty} e^{ -(\lambda_N + i\frac{y}{N} )k } b_{d}(k) d\hat{B}\big(N^{-1}y\big)\\
&=\frac{1}{N^{d-\frac{1}{2}}}\int_{\rr} \widehat{f_{N,m}}(y) \frac{\frac{iy}{N}}{e^{\frac{iy}{N}}-1} \sum_{k=0}^{\infty} e^{ -(\lambda_N + i\frac{y}{N} )k } b_{d}(k) d\hat{B}\big(N^{-1}y\big),
\end{split}
\end{equation}
where $\widehat{f_{N,m}}(y)=\sum_{j=0}^{m} f\Big(\frac{j}{N}\Big) \frac{1}{\sqrt{2\pi}} \int_{\rr}e^{i\omega y}\1_{(\frac{j}{N},\frac{j+1}{N})}(\omega)\ d\omega$ is the Fourier transform of ${f_{N,m}}$. We note 
\begin{equation}\label{upperbound for sum}
\sum_{k=0}^{\infty} e^{- (\lambda_N + \frac{iy}{N})k} b_{d}(k) < C \Big(\lambda_N + \frac{iy}{N}\Big)^{-d},
\end{equation}
for $d>-\frac{1}{2}$ and a constant $C$ by Lemma \ref{lem1}. Using \eqref{sumwelldefined} and \eqref{upperbound for sum}, we have
\begin{eqnarray}\label{upper bound expected second}
\mathbb{E}\Bigg| \frac{1}{ {N}^{d+ \frac{1}{2}} }\sum_{j=0}^{m}f\Big(\frac{j}{N}\Big) {X}_{d,\lambda_N}(j)\Bigg|^{2}
\!\!\!\!&=&\!\!\!\! \int_{\rr}\Big| \widehat{f_{N,m}}(y)\Big|^{2} \Bigg| \frac{\frac{iy}{N}}{e^{\frac{iy}{N}}-1} \Bigg|^{2}
\frac{1}{ {N}^{2d} }\ \!\!\Big| \sum_{k=0}^{\infty} e^{- (\lambda_N + \frac{iy}{N})k} b_{d}(k)   \Big|^{2}\ dy \nonumber \\
\!\!\!\! &\leq& \!\!\!\! \frac{\pi^2}{4}\, C \int_{\rr}\Big| \widehat{f_{N,m}}(y)\Big|^{2} \Big[ (N\lambda_N)^2 + y^2 \Big]^{-d}\ dy \label{inequality}\\
\!\!\!\! &=&\!\!\!\!  C' \|f_{N,m}\|^{2}_{{\mathcal{A}}_{3,N\lambda_N}}, \nonumber
\end{eqnarray}
where $C'$ is another constant. Now, for $m_2>m_1\geq 1$, we have
\begin{equation*}
\mathbb{E}\Bigg|\frac{1}{ {N}^{d+ \frac{1}{2}}   }\sum_{j=m_1+1}^{m_2}f\Big(\frac{j}{N}\Big) X_{d,\lambda_N}(j)\Bigg|^{2}\leq C'\|{f}^+_{N,m_2}-{f}^+_{N,m_1}\|^{2}_{{\mathcal{A}}_{3,\lambda_N}}\to 0
\end{equation*}
as $m_1,m_2\to \infty$ and this shows the series is well-defined. The following remark illustrates the inequality in~\eqref{inequality}.
\begin{rem}
In~\eqref{inequality} we used Lemma~\ref{lem1} and $\Big|\frac{\frac{iy}{N}}{e^{\frac{iy}{N}}-1} \Big|^{2} \leq \frac{\pi^2}{4}$ for $y \in [-N\pi,N\pi]$. This can be seen as follows
\begin{equation*}
\Bigg|\frac{\frac{iy}{N}}{e^{\frac{iy}{N}}-1} \Bigg|^{2} = \frac{\frac{y^2}{N^2}}{|\cos \frac{y}{N}+i\sin\frac{y}{N}-1|^2}
                                                           = \frac{\frac{y^2}{N^2}}{2(1-\cos \frac{y}{N})}  = \frac14 \frac{\frac{y^2}{N^2}}{\sin^2 \frac{y}{2N}}.
\end{equation*}
Then taking the limit yields
\begin{equation*}
\lim_{N\to\infty} \frac14 \frac{\frac{N^2\pi^2}{N^2}}{\sin^2 \frac{\pm N\pi}{2N}} = \frac14 \lim_{N\to\infty} \frac{\pi^2}{\sin^2 \pm\frac{\pi}{2}}= \frac{\pi^2}{4}.
\end{equation*}
\end{rem}

Next, we show that $U_{N}$ converges in distribution to $U$ as $N\to\infty$. By applying a similar proof to that of Theorem 3.15 in~\cite{MeerschaertsabzikarSPA}, it can be shown that the set of elementary functions are dense in ${\mathcal{A}}_{d,\lambda}$ and then there exists a sequence of elementary functions $f^{l}$ such that $\|f-f^{l}\|_{{\mathcal{A}}_{3}}\to 0$, as $l\to\infty$. Now, assume
\begin{equation}\label{eq:U defn}
U^{l}_{N}=\frac{1}{N^{ d+ \frac{1}{2}}}\sum_{j=-\infty}^{\infty} f^{l}\Big(\frac{j}{N}\Big) X_{d,\lambda_N}(j),\quad \ U^{l}=\Gamma^{-1}(d+1)\int_{\rr}f^{l}(u)\ dB^{I\!I}_{d,\lambda_*}(u).
\end{equation}
Observe that $U^{l}_{N}$ is well-defined, since $U^{l}_{N}$ has a finite number of terms and the elementary function $f^{l}$ is in ${\mathcal{A}}_{3}$. According to~\cite[Theorem 4.2.]{Bill}, the series $U_{N}$ converges in distribution to $U$ if
\begin{description}
  \item[Step 1] $U^{l}\longrightarrowD U$, as $l\to\infty$,
  \item[Step 2] for all $l\in\mathbb{N}$, $U^{l}_{N}\longrightarrowD U^{l}$, as $N\to\infty$,
  \item[Step 3] $\limsup_{l\to\infty}\limsup_{N\to\infty}\mathbb{E}\Big|U^{l}_{N}-U_{N}\Big|^{2}=0$.
\end{description}

Step $1$: The random variables $U^{l}$ and $U$ have normal distribution with mean zero and variances $\|f^{l}\|_{{\mathcal{A}}_{3,\lambda_N}}$ and $\|f\|_{{\mathcal{A}}_{3,\lambda_N}}$, respectively, since $f$ and $f^{l}$ are in ${{\mathcal{A}}_{3,\lambda_N}}$. Therefore $\mathbb{E}\Big|U^{l}-U\Big|^{2}=\|f^{l}-f\|_{{\mathcal{A}}_{3,\lambda_N}}\to 0$ as $l\to\infty$.

Step $2$: Note that $f^{l}$ is an elementary function and hence $U^{l}_{N}$, given by \eqref{eq:U defn}, can be written as
$U^{l}_{N}=\frac{1}{ N^{d+\frac{1}{2} } }\int_{\mathbb{R}} f^{l}(u) dS_{d,\lambda_N}(u)$.
Now, apply part (iii) of Theorem 4.3 in \cite{SabzikarSurgailisinvariance} to see that
$\frac{S_{d,\lambda_N}(u)}{ N^{d+\frac{1}{2} } }\limfdd \Gamma^{-1}(d+1) B^{I\!I}_{d,\lambda_*}(u)$, as $N\to\infty$, and this implies that $U^{l}_{N}\limfdd U^{l}$, as $N\to\infty$.

Step $3$: By a similar arguments of \eqref{upperbound for sum} and \eqref{upper bound expected second}, we have
\begin{equation}\label{eq:upperboundforeypectedU}
\begin{split}
\mathbb{E}\Big|U^{l}_{N}-U_{N}\Big|^{2}
&=\int_{\rr} \Big| \widehat{f^{l}_{N}}(y)-\widehat{f_N}(y)\Big|^{2}    \Big|\frac{\frac{iy}{N}}{e^{\frac{iy}{N}}-1} \Big|^{2}
\frac{1}{N^{2d}}\ \Big|1-e^{-(\lambda_N+\frac{iy}{N})}\Big|^{-2d}\ dy\\
&\leq C \int_{\rr} \Big| \widehat{f^{l}_{N}}(y)-\widehat{f_N}(y)\Big|^{2} \Big[ (N \lambda_N)^{2} + y^2 \Big]^{-d}\ dy,
\end{split}
\end{equation}
where $\widehat{f^{l}_{N}}(y)$ and $\widehat{f_N}(y)$ are the Fourier transforms of
\begin{equation*}
f^{l}_{N}(u):=\sum_{j=0}^{\infty}f^{l}\Big(\frac{j}{N}\Big){\bf 1}_{ \Big(\frac{j}{N},\frac{Nx-j+1}{N}\Big)}(u)
\end{equation*}
and $f_{N}:=\sum_{j=0}^{\infty}f\Big(\frac{j}{N}\Big) {\bf 1}_{ \big(\frac{j}{N},\frac{Nx-j+1}{N}\big) }(u)$ respectively. Note that $f^{l}$ is an elementary function and therefore $\widehat{f^{l}_{N}}$ converges to $\widehat{f^{l}}$ at every point and $\Big|\widehat{f^{l}_{N}}(\omega)-\widehat{f^{l}}(\omega)\Big|\leq \widehat{g^{l}}(\omega)$ uniformly in $N$, for some function $\widehat{g^{l}}(\omega)$ which is bounded by the minimum of $C_1$ and $C_2|\omega|^{-1}$ for all $\omega\in\rr$ (See Theorem 3.2. in \cite{PipirasTaqqu2} for more details). Let $\mu_{d,\lambda}(d\omega)= (\lambda^2+\omega^2)^{-d}\ d\omega$ be the measure on the real line for $d >-\frac{1}{2}$, then  $\widehat{g^{l}}(\omega)\in L^{2}(\mathbb{R},\mu_{d,\lambda})$. Now apply the Dominated Convergence Theorem to see that
\begin{equation}\label{eq:lastproof}
\|f^{l}_{N}-f^{l}\|^{2}_{{\mathcal{A}}_{3}}=\|\widehat{f^{l}_{N}}-\widehat{f^{l}}\|^{2}_{L^{2}(\rr,\mu_{d,\lambda_*})}\to 0,
\end{equation}
as $N\to\infty$. From \eqref{sumwelldefined} and \eqref{eq:lastproof}, we have
\begin{equation*}
\begin{split}
\mathbb{E}\Big|U^{l}_{N}-U_{N}\Big|^{2}&\leq C\|  \widehat{ f^l_N } - \widehat{f_N} \|^{2}_{{\mathcal{A}}_{3}}\\
&\leq C\Big[\| \widehat{ f^l_N }- \widehat{f^l }\|^{2}_{{\mathcal{A}}_{3}}+\|\widehat{f} - \widehat{f_N }\|^{2}_{{\mathcal{A}}_{3}}+\| \widehat{f^l }-\widehat{f}\|^{2}_{{\mathcal{A}}_{3}}\Big].
\end{split}
\end{equation*}
The first two terms tend to zero as $N\to\infty$ because of \eqref{eq:lastproof} and Condition A respectively, and the last term tends to zero as $l\to 0$ (see step 1) and this completes the proof of Step 3.  \hfill $\Box$

\noindent{\bf proof Theorem~\ref{thm:the third question}:}
We prove only part (c) and omit the proofs of parts (a) and (b) due to the similarity of proofs. We first show that
\begin{equation}\label{goal_convegence_to_TFBM_Nh}
\frac{1}{(Nh)^{d+\frac{1}{2}}}\sum_{j=1}^{N} K\Big(\frac{Nx-j}{Nh}\Big) X_{d,\lambda_N}(j)\limfdd \frac{1}{\Gamma(d+1)}\int_{0}^{2}K^{'}(1-t) B^{I\!I}_{d,\lambda_*}(t) dt,
\end{equation}
where $B^{I\!I}_{d,\lambda_*}(t)$ is TFBMII.  Starting from the l.h.s of \eqref{goal_convegence_to_TFBM_Nh}, using Riemann sums to integrals and integration by parts
\begin{equation}\label{step1_convegence_to_TFBM_Nh}
\begin{split}
\frac{1}{(Nh)^{d+\frac{1}{2}}} \sum_{j=1}^{N} K\Big(\frac{Nx-j}{Nh}\Big) X_{d,\lambda_N}(j) &=
\frac{1}{(Nh)^{d+\frac{1}{2}}} \int_{0}^{1} K\Big(\frac{x-y}{h}\Big) dS_{d,\lambda_N}(y)\\
& = \frac{1}{(Nh)^{d+\frac{1}{2}}} \int_{0}^{1} K^{'}\Big(\frac{x-y}{h}\Big) S_{d,\lambda_N}(y) dy\\
&= \frac{1}{(Nh)^{d+\frac{1}{2}}} \int_{-1}^{1} K^{'}(u) S_{d,\lambda_N}(x-hu) du\\
& = \frac{1}{(Nh)^{d+\frac{1}{2}}} \int_{-1}^{1} K^{'}(u)  \sum_{j=\lfloor N(x-h)\rfloor}^{\lfloor N(x-hu)\rfloor} X_{d,\lambda_N}(j) du,
\end{split}
\end{equation}
where we used
\begin{equation*}
S_{d,\lambda_N}(x-hu) =\sum_{j=1}^{\lfloor N(x-h)\rfloor -1}X_{d,\lambda_N}(j)+ \sum_{j=\lfloor N(x-h)\rfloor}^{\lfloor N(x-hu)\rfloor} X_{d,\lambda_N}(j)
\end{equation*}
and the assumptions on the kernel function $K$ to see that
\begin{equation*}
\frac{\sum_{j=1}^{\lfloor N(x-h)\rfloor -1}X_{d,\lambda_N}(j)}{(Nh)^{d+\frac{1}{2}}} \int_{-1}^{1} K^{'}(u) du=0.
\end{equation*}
Next, by stationarity of $X_{d,\lambda_N}(j)$ and a change of variable we have
\begin{equation}\label{step2_convegence_to_TFBM_Nh}
\begin{split}
\sum_{j=\lfloor N(x-h)\rfloor}^{\lfloor N(x-hu)\rfloor} X_{d,\lambda_N}(j) &= \sum_{j=1}^{ \lfloor N(x-hu)\rfloor - \lfloor N(x-h)\rfloor +1} X_{d,\lambda_N}(j+ \lfloor N(x-h)\rfloor -1)\\
&\eqfdd \sum_{j=1}^{ \lfloor N(x-hu)\rfloor - \lfloor N(x-h)\rfloor +1} X_{d,\lambda_N}(j)= \sum_{j=1}^{ l_{x}(u) } X_{d,\lambda_N}(j),
  \end{split}
  \end{equation}
where $l_{x}(u)= \lfloor N(x-hu)\rfloor - \lfloor N(x-h)\rfloor +1$. Consequently, from \eqref{step1_convegence_to_TFBM_Nh} and \eqref{step2_convegence_to_TFBM_Nh}, we get
\begin{eqnarray}\label{step3_convegence_to_TFBM_Nh}
\frac{1}{(Nh)^{d+\frac{1}{2}}} \!\!\int_{-1}^{1} \!\!\!K^{'}(u)\!\!\!\!\!\!\!\sum_{j=\lfloor N(x-h)\rfloor}^{\lfloor N(x-hu)\rfloor}\!\!\!\!\!\!\! X_{d,\lambda_N}(j)\ du \!\!\!\!&\eqfdd& \!\!\!\!\!\!\!\! \frac{1}{(Nh)^{d+\frac{1}{2}}}\!\! \int_{-1}^{1}\!\!\! K^{'}(u)\sum_{j=1}^{ l_{x}(u) } X_{d,\lambda_N}(j)\ du \nonumber \\
\!\!\!\!\!&=&\!\!\!\!\!\!\!\!\frac{1}{(Nh)^{d+\frac{1}{2}}} \!\!\int_{0}^{2}\!\!\!\! K^{'}(1-t)\!\!\!\sum_{j=1}^{ \lfloor Nht\rfloor } \!\!\!X_{d,\lambda_N}(j) dt + o_{p}(1)
\end{eqnarray}
According to~\cite[Theorem 4.3]{SabzikarSurgailisinvariance} we have
  \begin{equation}\label{step4_convegence_to_TFBM_Nh}
\frac{1}{(Nh)^{d+\frac{1}{2}}}  \sum_{j=1}^{ \lfloor Nh  t\rfloor } X_{d,\lambda_N}(j)\limfdd \frac{1}{\Gamma(d+1)}B^{I\!I}_{d,\lambda_*}(t)
 \end{equation}
in $D[0,2]$ provided $Nh\to\infty$. Now, the desired result \eqref{goal_convegence_to_TFBM_Nh} follows from \eqref{step3_convegence_to_TFBM_Nh}, \eqref{step4_convegence_to_TFBM_Nh}, and the continuous mapping theorem. \hfill $\Box$

{\bf proof of Theorem \ref{flt}:} The proof of this theorem follows by Proposition \ref{lem:g squar integrable} and Theorem \ref{thm:the third question} and hence we omit the details. \hfill $\Box$

{\bf proof of Theorem \ref{asymptoticsNW}:}

For brevity, we restrict the proof of the theorem to $k=2$. Moreover, we just prove part (c) of the theorem since the proofs of the other two cases are similar. We first note that for each $0<x_i<1$,
\begin{equation}\label{Ahatdefn}
\widehat{A}_{N,i}= (Nh)^{\frac{1}{2}-d}[\hat{m}(x_i)-\mathbb{E}\hat{m}(x_i)]=\frac{1}{(Nh)^{d+\frac{1}{2}}}\sum_{1\leq s \leq N}K\Big(\frac{Nx_i -s}{Nh}\Big)X_{d,\lambda_N}(s).
\end{equation}
Let $j_i$ be integers such that $|Nx_i-j_i|\leq 1$ for $i=1,2$ and define
\begin{equation}\label{Atildedefn}
\widetilde{A}_{Ni}= \frac{1}{(Nh)^{d+\frac{1}{2}}}\sum_{s=|j_i|-\lfloor Nh \rfloor}^{|j_i|+\lfloor Nh \rfloor } K\Big(\frac{j_i -s}{Nh}\Big)X_{d,\lambda_N}(s)
\end{equation}
for $i=1, 2$. Since the kernel function $K$ vanishes in $\rr\setminus[-1,1]$ and $|K^{'}(x)|\leq C$ for all $x\in [-1,1]$, it follows that
\begin{equation}\label{Ahattildediff}
\widehat{A}_{N,i}-\widetilde{A}_{N,i}= o_{p}(1)
\end{equation}
for $i=1,2$. By a change of variable $s=\nu+j_1 -Nh$ and $s=\nu+j_2 -Nh-\lfloor N\delta\rfloor$, with $\delta= x_2-x_1$, for $\widetilde{A}_{N,1}$ and $\widetilde{A}_{N,2}$ respectively, use the fact that $X_{d,\lambda_N}$ is stationary, and $|K^{'}(x)|\leq C$ to see that
\begin{equation}\label{Aasttildediff}
\Big( \widetilde{A}_{N,1}, \widetilde{A}_{N,2} \Big)\eqfdd \Big( A^{\ast}_{N,1}, A^{\ast}_{N,2} \Big)+ o_{p}(1),
\end{equation}
where
\begin{equation}\label{AhatdefnN1}
A^{\ast}_{N,1} = \sum_{\nu=1}^{2\lfloor Nh\rfloor} K\Big(\frac{\nu}{Nh}-1\Big) X_{d,\lambda_N}(\nu)
\end{equation}
and
\begin{equation}\label{AhatdefnN2}
A^{\ast}_{N,2} = \sum_{\nu=\lfloor N\delta \rfloor}^{ \lfloor N\delta \rfloor + 2\lfloor Nh\rfloor} K\Big(\frac{\nu- \lfloor N\delta \rfloor}{Nh}-1\Big) X_{d,\lambda_N}(\nu).
\end{equation}
We use the partial sums $A^{\ast}_{N,i}$, for $i=1,2$, to establish the functional limit theorems. Let $\{ K_m\}$ be a sequence of elementary functions
such that $K_m\to K$ in $L^{2}$ as $m\to\infty$. Define $A^{\ast}_{m,Ni}$ be as \eqref{AhatdefnN1} and \eqref{AhatdefnN2} with $K_m(x)=\sum_{i=1}^{m} a_i {\bf 1}( t_{i-1},t_{i} )(x)$, where $a_i$ are some constants and $-1\leq t_i \leq 1$ for $i=0, \ldots, m$. We can rewrite $A^{\ast}_{m,N,i}$ as
\begin{equation}\label{AhatdefnNewrep}
A^{\ast}_{m,N,i} = 2^{d+\frac{1}{2}} \int_{0}^{2} K_{m}(u-1) dS^{\ast}_{Ni}(u) + o_{p}(1),
\end{equation}
where
\begin{equation}\label{SastfnN1}
S^{\ast}_{N,1}(s) = \frac{1}{(2Nh)^{d+\frac{1}{2}}}\sum_{t=1}^{\lfloor\lfloor Nh\rfloor s\rfloor}  X_{d,\lambda_N}(t)
\end{equation}
and
\begin{equation}\label{SastfnN2}
S^{\ast}_{N,2}(s) = \frac{1}{(2Nh)^{d+\frac{1}{2}}}\sum_{t=1}^{\lfloor\lfloor Nh\rfloor s\rfloor}  X_{d,\lambda_N}(t+ \lfloor N\delta \rfloor).
\end{equation}
Using~\cite[Theorem 4.3]{SabzikarSurgailisinvariance} and the continuous mapping theorem yields
\begin{equation}\label{AhatdefnNewrep}
A^{\ast}_{m,N,i} \limfdd A^{\ast}_{m} = \int_{0}^{2} K_{m}(u-1) dB^{I\!I}_{d,\lambda_*}(u),
\end{equation}
as $N\to\infty$ and hence
\begin{equation}\label{varianceAast}
\begin{split}
\sigma^{2}_{ii} &= \int_{0}^{2} \int_{0}^{2} K_{m}(u-1) K_{m}(v-1) {\rm Cov}\Big( B^{I\!I}_{d,\lambda_*}(u), B^{I\!I}_{d,\lambda_*}(v) \Big) du\ dv\\
&= \int_{0}^{2} \int_{0}^{2} K_{m}(u-1) K_{m}(v-1) |u-v|^{d-\frac{1}{2}} K_{d-\frac{1}{2}}(\lambda_* |u-v|) du\ dv.
\end{split}
\end{equation}
Next, we need to show that $A^{\ast}_{m,N,1}$ and $A^{\ast}_{m,N,2}$ are asymptotically independent (i.e. $\sigma^{2}_{12} = \sigma^{2}_{21} =0$ ). Observe that
\begin{equation}
A^{\ast}_{m,N,1} = 2^{d+\frac{1}{2}} \sum_{j=1}^{m} a_{j} [ S^{\ast}_{N1}(t_{j}) - S^{\ast}_{N1}(t_{j-1})] + o_{P}(1)
\end{equation}
and
\begin{equation}\label{SastfnN1diff}
S^{\ast}_{N,1}(t_j) - S^{\ast}_{N,1}(t_{j-1}) = \frac{1}{(2Nh)^{d+\frac{1}{2}}}\sum_{s=\lfloor\lfloor Nh \rfloor t_{j-1}\rfloor}^{\lfloor\lfloor Nh\rfloor t_j\rfloor}  X_{d,\lambda_N}(s)\\
=\sum_{p=-\infty}^{\infty} d_{pN}\zeta(p),
\end{equation}
where
\begin{equation}\label{dpNdefn}
d_{pN}=\frac{1}{(2Nh)^{d+\frac{1}{2}}} \sum_{t= \lfloor\lfloor Nh\rfloor t_{j-1}\rfloor}^{\lfloor\lfloor Nh\rfloor t_j\rfloor} b_{d}(p-t)e^{-\lambda_N (p-t)}.
\end{equation}
Since $b_{d}(j)e^{-\lambda_N j}\sim C j^{d-1} e^{-\lambda_N j}$ for large lag $j$, see \eqref{bdk}, then for $p> {\lfloor\lfloor Nh\rfloor t_j\rfloor}$, we have
\begin{equation}\label{dpNupperbound}
|d_{pN}| \leq C |Nh|^{d+\frac{1}{2}} Nh (p-{\lfloor\lfloor Nh\rfloor t_j\rfloor})^{d-1} e^{-\lambda_N (p-{\lfloor\lfloor Nh\rfloor t_j\rfloor})},
\end{equation}
where $C$ is a constant. Therefore, we get
\begin{equation}\label{dpNsquaresum}
\lim_{N \to\infty}\sum_{|p|>M} d^2_{pN}= 0,
\end{equation}
since $h\log(Nh) \to 0$ and $M=Nh\log(Nh)$. Consequently,
\begin{equation}\label{SastN1square}
\lim_{N\to\infty}\mathbb{E}\Bigg[ S^{\ast}_{N1}(t_j) - S^{\ast}_{N1}(t_{j-1}) - \sum_{|p|\leq M} d_{pN} \zeta(p)      \Bigg]^{2} = 0
\end{equation}
and by a similar argument
\begin{equation}\label{SastN2square}
\lim_{N\to\infty}\mathbb{E}\Bigg[ S^{\ast}_{N2}(t_j) - S^{\ast}_{N2}(t_{j-1}) - \sum_{|p|\leq M} d_{pN} \zeta(p+\lfloor N\delta\rfloor)     \Bigg]^{2} = 0
\end{equation}
From \eqref{SastN1square}, \eqref{SastN2square}, and $ \lfloor N\delta \rfloor - 2N\to\infty$, we conclude that
$ S^{\ast}_{N,1}(t_j) - S^{\ast}_{N,1}(t_{j-1})$ and $ S^{\ast}_{N,2}(t_{j^{'}}) - S^{\ast}_{N,2}(t_{{ j^{'} } -1})$ are asymptotically independent for all
$j, j^{'}$ and this implies that $A^{\ast}_{m,N,1}$ and $A^{\ast}_{m,N,2}$ are asymptotically independent. Thus
\begin{equation}
\Big(A^{\ast}_{m,N,1}, A^{\ast}_{m,N,2}\Big)\limfdd N_{2}\Big( 0, \Sigma \Big),
\end{equation}
where $\sigma^{2}_{ii}$ is given by \eqref{varianceAast} and $\sigma^{2}_{12}=\sigma^{2}_{21}=0$. Since
\begin{equation}\label{varianceA}
\lim_{m\to\infty}\sigma^{2}_{ii} = \int_{0}^{2} \int_{0}^{2} K(u-1) K(v-1) |u-v|^{d-\frac{1}{2}} K_{d-\frac{1}{2}}(\lambda_* |u-v|) du\ dv
\end{equation}
and by~\cite[Theorem 2]{Deo2}, one obtains
\begin{equation}\label{DeoTheorem2}
\lim_{m\to\infty}\lim_{N\to\infty} {\rm Var}(A^{\ast}_{m,N,1}- A^{\ast}_{N,1})=0,\ {\rm for}\ i=1, 2,
\end{equation}
then the desired results follows by \eqref{Ahattildediff}, \eqref{Aasttildediff}, \eqref{varianceA}, and \eqref{DeoTheorem2}. $\hfill \Box$

{\bf proof of Theorem \ref{Theorem1}:} The proofs of part (a) and part (b) are similar and we just proof part (b) for $\lambda_*\in (0,\infty)$. Using the triangle inequality yields
\begin{eqnarray*}
&& \mathbb{P}\Big( N^{\frac{1}{2} - d}\big\|\hat{\theta}-\theta- (\M_{N+}^T\M_{N+})^{-1} \M_{N+}^T\e_{N}\big\| >\Delta \Big) \\ &\leq& \mathbb{P}\Big( N^{\frac{1}{2} - d}\Big\|\hat{\theta}-\theta\Big\| >\frac{\Delta}{2} \Big)                                                      +\mathbb{P}\Big( N^{\frac{1}{2} - d}\big\|(\M_{N+}^T\M_{N+})^{-1} \M_{N+}^T\e_{N}\big\| >\frac{\Delta}{2} \Big).
\end{eqnarray*}
For the first term, we have by~\cite[p. 95 Theorem 5.1.9]{Weiershauser2012}, ~\eqref{semilongmemory} and Markov's inequality,
\begin{equation*}
\mathbb{P}\big( N^{\frac{1}{2} - d}\big\|\hat{\theta}-\theta\big\| >\frac{\Delta}{2} \Big) \leq \frac{4\E\big\|\hat{\theta}-\theta\big\|^2}{\Delta^2 N^{1-2d}} = \frac{O[\min(N^{-\beta}N^{2d-1},N^{4d-2})]}{\Delta^2}.
\end{equation*}
Since $\Delta$ is arbitrary, $\lim_{N\to\infty}\mathbb{P}\big( N^{\frac{1}{2} - d}\big\|\hat{\theta}-\theta\big\| >\frac{\Delta}{2} \Big) = 0$ for $d < 1/2$. For the second term and again by Markov's inequality we have
\begin{eqnarray}
\mathbb{P}\Big(\! N^{\frac{1}{2} - d}\big\|(\M_{N+}^T\M_{N+})^{-1} \M_{N+}^T\e_{N}\big\| \!>\!\frac{\Delta}{2} \Big) \!\!\!&\leq& \!\!\! \frac{4\E\big\|(\M_{N+}^T\M_{N+})^{-1} \M_{N+}^T\e_{N}\!\big\|^2}{\Delta^2 N^{1-2d}}.
\label{proof53}
\end{eqnarray}
Let $\Omega = (\M_{N+}^T\M_{N+})^{-1} \M_{N+}^T$, then $\E\|\Omega \e_N\|^2 = \tr(\Omega \Sigma_{\e_N\e_N}) + \{\E(\e_N)\}^T\Sigma_{\e_N\e_N}\E(\e_N)$ where $\tr(A)$ and $\Sigma_{\e_N\e_N}$ denote the trace of the matrix $A$ and the variance-covariance matrix of $\e_N$ respectively. Since $\{ X_{d,\lambda}(t) \}_{t\in\Z}$ is a tempered mean zero linear process we have $\E\|\Omega \e_N\|^2 = \tr(\Omega \Sigma_{\e_N\e_N})$. Further, the variance-covariance matrix of $\e_N$ is finite, see~\eqref{covZ}. Consequently, the numerator of~\eqref{proof53} is finite since $\M_{N+}^T\M_{N+}$ is full rank for $N\to\infty$ and hence the second term goes to zero for $d<1/2$. \hfill $\Box$

{\it proof of Theorem \ref{Theorem4}:}
The proofs of part (a) and part (b) are similar and hence we just give the proof for part (b). We first note that $\int_\rr \mu_{(i+)}(u) dB^{I\!I}_{d,\lambda_*}(u) = \int_{\rr} {\mathbb I}^{d,\lambda_*}_{-}\mu_{(i+)}(u) dB(u) $. Observe that $\mu_{(i+)}\in L^{p}(\rr)$
for $p\geq 1$ and hence ${\mathbb I}^{d,\lambda_*}_{-}\mu_{(i+)}\in L^{p}$. In particular, let $p=2$ and apply the Ito-isometry to conclude that
$\int_{\rr} {\mathbb I}^{d,\lambda_*}_{-}\mu_{(i+)}(u) dB(u)$ is well-defined. Because of Theorem~\ref{Theorem1}, we need to show that
\begin{equation}\label{th1step1}
N^{\frac{1}{2} - d} (\M_{N+}^T\M_{N+})^{-1} \M_{N+}^T \e_N \limfdd {\bf \Lambda} \Big[ \int_\rr \mu_{(i+)}(u)\ dB^{I\!I}_{d,\lambda_*}(u) \Big]_{i= 1, \ldots, p+1}
\end{equation}
as $N\to\infty$. Observe that $N (\M_{N+}^T\M_{N+})^{-1} \to \Lambda$ as $N\to\infty$. Therefore the RHS of \eqref{th1step1} follows if we show
\begin{equation*}
\frac{1}{N^{ d+ \frac{1}{2} } } \M_{N+}^T \e_N \limfdd \Big[\int_{\rr} \mu_{(i+)}(s)\ dB^{I\!I}_{d,\lambda_*}(s)\Big]_{i=1,\ldots, p+1}.
\end{equation*}
But this is equivalent to show that
\begin{equation}
\frac{1}{N^{ d+ \frac{1}{2} } } {\langle \alpha, \M_{N+}^T \e_N \rangle} \to {\Bigl\langle \alpha, \Big[ \int_\rr \mu_{(i+)}(u)\ dB^{I\!I}_{d,\lambda_*}(u) \Big]_{i= 1, \ldots, p+1} \Bigr\rangle}, \alpha\in\rr^{p+1}.
\end{equation}
Note that
\begin{equation}
{\langle \alpha, \M_{N+}^T \e_N \rangle} = \sum_{i=1}^{p+1}\sum_{j=1}^{n} \alpha_i \mu_{(i+)}\big(\frac{j}{N}\big) X_{d,\lambda_N}(j)
\end{equation}
and by Lemma~\ref{thm:similar the third question} $$N^{-(d+1/2)}\sum_{i=1}^{p+1}\sum_{j=1}^{N} \alpha_i \mu_{(i+)}(\frac{j}{N}) X_{d,\lambda_N}(j)\limfdd \int_{\rr} m_{\alpha}(u) dB^{I\!I}_{d,\lambda_*}(u)$$
where $m_{\alpha}(u) := \sum_{i=1}^{p+1} \alpha_i \mu_{(i+)}(u)  $ and this completes the proof.   \hfill $\Box$

{ \bf proof of Theorem \ref{limitdistributiontheorem}: }
We only proof part (c) since the other parts follow a similar procedure. Let $\Xi$ be the random vector
\begin{equation}
\Xi = \Big[\int_{\rr} \mu_{(i+)}(s) dB^{I\!I}_{d,\lambda_*}(s)\Big]_{i=1,\ldots, p+1}.
\end{equation}
Then we can write
\begin{equation}
\int_\rr \mu_{(i+)}(s) dB^{I\!I}_{d,\lambda_*}(s) = \int_\rr \Big( {\mathbb I}^{d,\lambda_*}_{-}\mu_{(i+)} \big)(s) dB(s)
\end{equation}
for $i=1,\ldots, p+1$. We observe that $\int_\rr \big( {\mathbb I}^{d,\lambda_*}_{-}\mu_{(i+)} \big)(s) dB(s)$ is a Gaussian stochastic process with mean zero and
finite variance $\int_\rr \big| {\mathbb I}^{d,\lambda_*}_{-}\mu_{(i+)}(s) \big|^2 \ ds$. Using the Ito-isometry for the Wiener integrals, one can see $\Xi$ has the covariance matrix
\begin{equation}
\Sigma_{0}=\Bigg [\int_{\mathbb{R}} \Big( {\mathbb I}^{d,\lambda_*}_{-} \mu_{(i+)} \Big)(s)
 \Big( {\mathbb I}^{d,\lambda_*}_{-} \mu_{(k+)} \Big)(s)\ ds\Bigg]_{i,k=1,\ldots,p+1}
\end{equation}
and consequently $\Lambda \Xi$ has normal distribution with covariance matrix $\Lambda \Sigma_0 \Lambda$ and this completes the proof of the first part. Next, we have
\begin{equation}\label{covariance}
\begin{split}
\int_{\mathbb{R}} \Big( {\mathbb I}^{d,\lambda_*}_{-} \mu_{(i+)} \Big)(s)
 \Big( {\mathbb I}^{d,\lambda_*}_{-} \mu_{(k+)} \Big)(s)\ ds&=
\int_{\rr} \mathcal{F}[{\mathbb I}^{d,\lambda_*}_{-} \mu_{(i+)}](\omega)
\overline{\mathcal{F}[{\mathbb I}^{d,\lambda_*}_{-} \mu_{(i+)}](\omega)} \ d\omega\\
&= \int_\rr  \widehat{\mu_{(i+)}}(\omega) \overline{\widehat{\mu_{(k+)}}(\omega)}  (\lambda_*^2 + \omega^2)^{-d}d\omega\\
&=\int_\rr \int_\rr \mu_{(i+)}(t) \mu_{(k+)}(s) \int_{\rr} e^{i\omega(t-s)} (\lambda_*^2 + \omega^2)^{-d} \ d\omega ds\ dt\\
&=2\int_\rr \int_\rr \mu_{(i+)}(t) \mu_{(k+)}(s) \int_{0}^{\infty} \cos(\omega(t-s)) (\lambda_*^2 + \omega^2)^{-d} d\omega\ ds\ dt\\
&=C \int_\rr \int_\rr \mu_{(i+)}(t) \mu_{(k+)}(s) |t-s|^{d-\frac{1}{2}} K_{d-\frac{1}{2}}(\lambda_*|t-s|) ds \ dt,
\end{split}
\end{equation}
where $C = \frac{2}{\Gamma(d) \sqrt{\pi} (2\lambda)^{d-\frac{1}{2}}}$ and we used
\begin{equation}
\int_{0}^{\infty} \frac{\cos(\omega x)}{(\lambda^2+x^2)^{\nu+\frac{1}{2}}   } \ dx = \frac{\sqrt{\pi}}{\Gamma(\nu+\frac{1}{2})} \Big(  \frac {|x|}{2\lambda} \Big)^{\nu} K_{\nu}(\lambda|x|)
\end{equation}
for $\nu > -\frac{1}{2}$ and $\lambda>0$ in \eqref{covariance}. This completes the proof of the second part and Theorem.  \hfill $\Box$

\end{document}